\documentclass[11pt,letterpaper]{article}
\usepackage{amsmath}
\usepackage{amsfonts}
\usepackage{amssymb}
\usepackage{fullpage}
\usepackage[top=1.2in, bottom=1.4in, left=1.2in, right=1.2in]{geometry} 
\usepackage{setspace}
\usepackage[pdftex,breaklinks]{hyperref}

\newcommand{\beq}{\begin{equation}}
\newcommand{\eeq}{\end{equation}}
\newcommand{\bit}{\begin{itemize}}
\newcommand{\eit}{\end{itemize}}
\newcommand{\ben}{\begin{enumerate}}
\newcommand{\een}{\end{enumerate}}
\newcommand{\tm}{\textrm}

\hypersetup{
    unicode=false,          
    pdftoolbar=true,        
    pdfmenubar=true,        
    pdffitwindow=false,     
    pdfstartview={Fit},     
		pdfpagelayout={SinglePage},		
    pdftitle={},    							
    pdfauthor={Anshul Agarwal},		
    pdfsubject={Documentation},		
    pdfcreator={Anshul Agarwal},	
    pdfproducer={Anshul Agarwal},	
    pdfkeywords={}, 							
    pdfnewwindow=true,      
    colorlinks=true,        
    linkcolor=blue,         
    citecolor=blue,         
    filecolor=blue,         
    urlcolor=blue           
}

\title{A Novel MINLP Reformulation for Nonlinear Generalized Disjunctive Programming (GDP) Problems}

\author{Anshul Agarwal\thanks{The Dow Chemical Company, 2301 N. Brazosport Blvd., Freeport, TX 77541, United States}}

\begin{document}

\maketitle
\date{}

\begin{abstract}

In optimization problems, often equations and inequalities are represented using if-else (implication) construct which is known to be equivalent to a disjunction. Such statements are modeled and incorporated in an optimization problem using Generalized Disjunctive Programming (GDP).  GDP provides a systematic methodology to model optimization problems involving logic disjunctions, logic propositions, and algebraic equations.  In order to take advantage of the existing MINLP solvers, GDP problems can be reformulated as the standard MINLP problems.  In this work we propose a novel reformulation methodology for general GDP problems with nonlinear equality and inequality constraints.  The proposed methodology provides an exact reformulation, maintains feasibility and convexity of the constraints, and, most importantly, does not require choosing a tolerance level and a Big-M parameter. We also demonstrate how the new reformulation approach can be used to convert the logic proposition represented using if-else (implication) construct into equations in the standard MINLP format. The conversion methodology is extended for variations of implication constructs that include implicit else blocks, sequential implication logic, multiple testing conditions, and nested implication blocks.  The proposed approach is utilized to model physical and mechanical properties in a mathematical optimization tool that solves an MINLP problem to design commercial products.

\end{abstract}

\section{Introduction}

\subsection{Generalized Disjunctive Programming (GDP)}

Formulating a nonlinear or mixed-integer nonlinear mathematical optimization problem often requires incorporating certain logical conditions together with other discrete/continuous decisions. These logical conditions are statements about equations and inequalities that involve operations such as ``and'' (conjunction), ``or'' (disjunction), and ``complement of'' (negation). In particular, often equations and inequalities are modeled using ``if$\ldots$then$\ldots$else'' (implication) construct which is known to be equivalent to a disjunction. For example, we can have a few equations in an optimization problem defined in the following manner:

\mbox{ }\\
\hspace*{1cm}if $T \geq \alpha$ then\\
\hspace*{2cm}if $g_{1}(T) \leq 0$ or $g_{2}(T) \leq 0$ then\\
\hspace*{3cm}$P = f_{1}(x)$\\
\hspace*{2cm}else\\
\hspace*{3cm}$P = f_{2}(x)$\\
\hspace*{2cm}end\\
\hspace*{1cm}else\\
\hspace*{2cm}$P = C$\\
\hspace*{1cm}end\\
\mbox{ }\\
Here what property model for mechanical property $P$ is used for optimization depends on the range in which the other mechanical property $T$ falls. Such statements can be modeled and incorporated in an optimization problem using Generalized Disjunctive Programming (GDP).

Generalized Disjunctive Programming, an extension of the disjunctive programming (DP)  developed by Balas \cite{balas1979, balas1985}, provides a systematic methodology to model optimization problems involving logic disjunctions, logic propositions, and algebraic equations \cite{raman1993, raman1994, beau1991}. A GDP problem can be regarded as a mixed-integer linear/nonlinear program with disjunctive constraints. GDP based representations help retain and exploit the inherent logic structure of problems that, as a result, reduce the combinatorics and improve relaxations and bounds of the global optimum, especially in nonconvex problems \cite{gross2012}.

In GDP problems, in general, disjunctions are represented as follows:

\beq
\label{eq:GDPgeneral}
	\begin{array}{l}
		\bigvee\limits_{j \in D_{k}}
		\left [
			\begin{array}{c}
				Y_{jk} \\
				h_{jk}(x) \leq 0
			\end{array}
		\right ], \quad k \in K \\
		\Omega(Y) = True \\
		x^{L} \leq x \leq x^{U}, \;\;\; Y_{jk} \in \{True,False\} 
	\end{array}
\eeq
Here we have a set of $K$ disjunctions (logical conditions). Each disjunction comprises a number of terms $j \in D_{k}$. In each term there is a Boolean variable $Y_{jk}$ and a set of inequalities $h_{jk} \leq 0$. If $Y_{jk}$ is true, inequalities $h_{jk} \leq 0$ are enforced, otherwise ignored. Also, $\Omega(Y) = True$ are logic propositions for the Boolean variables. All terms in a disjunction are connected by the logical ``or'' operator ($\bigvee$) that is inclusive or exclusive depeding on the constraints written on the Boolean variables. The constraints $h_{jk}$ can be linear, convex nonlinear, or nonconvex nonlinear in nature. We note that the vector $x$ includes variables that depend on which term in a disjunction is true as well as those that depend on equations and inequalities located outside of the disjunctions.

To illustrate the concept, the aforementioned implication construct for mechanical property $P$ can be represented using nested disjunctions in the following manner \cite{vecchietti2000}.

\begin{equation*}
	\begin{array}{l}
		\left [
			\begin{array}{c}
				Y_{1} \\
				T \geq \alpha \\
				\left [
					\begin{array}{c}
						Z_{1} \\
						\left [
							\begin{array}{c}
								W_{1} \\
								g_{1}(T) \leq 0 \\
								P = f_{1}(x)
							\end{array}
						\right ] \bigvee
						\left [
							\begin{array}{c}
								W_{2} \\
								g_{2}(T) \leq 0 \\
								P = f_{1}(x)
							\end{array}
						\right ]
					\end{array}
				\right ] \bigvee
				\left [
					\begin{array}{c}
						Z_{2} \\
						P = f_{2}(x)
					\end{array}
				\right ]
			\end{array}
		\right ] \bigvee
		\left [
			\begin{array}{c}
				Y_{2} \\
				P = C
			\end{array}
		\right ] \\
		Y_{1} = True \;\;\tm{or}\;\; Y_{2} = True \;\;\tm{not both} \\
		Z_{1} = True \;\;\tm{or}\;\; Z_{2} = True \;\;\tm{not both} \\
		W_{1} = True \;\;\tm{or}\;\; W_{2} = True \;\;\tm{or both} \\
		Y_{1} = True \Rightarrow Z_{1} = True \;\;\tm{or}\;\; Z_{2} = True \;\;\tm{not both} \\
		Z_{1} = True \Rightarrow W_{1} = True \;\;\tm{or}\;\; W_{2} = True \;\;\tm{or both}
	\end{array}
\end{equation*}
These nested disjunctions can eventually be converted into a set of three standard disjunctions as below \cite{vecchietti2000}. We discuss such a conversion in more detail in Section 3.

\begin{equation*}
	\begin{array}{l}
		\left [
			\begin{array}{c}
				Y_{1} \\
				T \geq \alpha \\
				Z_{1} \vee Z_{2}
			\end{array}
		\right ] \bigvee
		\left [
			\begin{array}{c}
				Y_{2} \\
				T < \alpha \\
				P = C
			\end{array}
		\right ] \\
		\left [
			\begin{array}{c}
				Z_{1} \\
				P = f_{1}(x) \\
				W_{1} \vee W_{2}
			\end{array}
		\right ] \bigvee
		\left [
			\begin{array}{c}
				Z_{2} \\
				P = f_{2}(x) \\
				g_{1}(T) > 0 \\
				g_{2}(T) > 0
			\end{array}
		\right ] \\
		\left [
			\begin{array}{c}
				W_{1} \\
				g_{1}(T) \leq 0
			\end{array}
		\right ] \bigvee
		\left [
			\begin{array}{c}
				W_{2} \\
				g_{2}(T) \leq 0
			\end{array}
		\right ] \\
	\end{array}
\end{equation*}

GDP problems can be solved using a specialized disjunctive Branch and Bound method \cite{lee2000}, or using the Logic-Based Outer Approximation method \cite{turkay1996} which efficiently exploits the logic structure of a GDP problem. On the other hand, in order to take advantage of the existing MINLP solvers, GDP problems can also be reformulated as the standard MINLP problems using either the Big-M reformulation \cite{nemhau1988, raman1991, will1999} or the Convex Hull reformulation \cite{balas1985, raman1994, gross2002, gross2003}. In this work, we utilize the approach of reformulating GDP problems into MINLP problems.

\subsection{GDP Convex Hull Reformulation and its Limitations}

In order to reformulate a general nonlinear nonconvex GDP problem into an MINLP problem, Lee and Grossmann \cite{lee2000} proposed a convex hull reformulation methodology. To illustrate their method, consider the following simple disjunction with two terms for two nonlinear equations:

\begin{equation*}
	\begin{array}{l}
		\left [ 
			\begin{array}{c}
				Y_{1} \\
				h_{1}(x) = 0
			\end{array}
		\right ] \vee 
		\left [ 
			\begin{array}{c}
				Y_{2} \\
				h_{2}(x) = 0
			\end{array}
		\right ] \\
		0 \leq x \leq x^{U}
	\end{array}
\end{equation*}
Here, if $Y_{1}$ is true, equation $h_{1}(x)$  will be satisfied, otherwise equation $h_{2}(x)$. According to the methodology proposed by Lee and Grossmann, a convex hull reformulation of this disjunction will be written as:

\begin{equation*}
	\begin{array}{l}
		x = \nu_{1} + \nu_{2} \\
		0 \leq \nu_{j} \leq \lambda_{j}x^{U}, \;\; j \in \{1,2\} \\
		\lambda_{1} + \lambda_{2} = 1 \\
		\lambda_{j}h_{j}(\nu_{j}/\lambda_{j}) = 0, \;\; j \in \{1,2\} \\
		x, \nu_{j} \geq 0, \;\; \lambda_{j} \in \{0,1\}, \;\; j \in \{1,2\}
	\end{array}
\end{equation*}
Here $\nu_{j}$ are disaggregated variables that are assigned to each term of the disjunction, and $x^{U}$ serve as their upper bounds. The binary variables $\lambda_{j}$ are the weight factors that determine the feasibility of the disjunctive term. When $\lambda_{j}$ is 1, the $j$'th constraint in the disjunction is enforced and other constraints are ignored. Lee and Grossmann also showed that the constraint $\lambda_{j}h_{j}(\nu_{j}/\lambda_{j})$ is convex if $h_{j}(x)$ is convex. Note that $h_{j}(x)$ can be an inequality as well.

In general, a set of nonlinear nonconvex disjunctions as shown in (\ref{eq:GDPgeneral}) can be reformulated into the following MINLP statements using the convex hull approach:
\beq
\label{eq:GDPCHreform}
	\begin{array}{l}
		\displaystyle{x = \sum_{j \in D_{k}} \nu_{jk}}, \;\; k \in K \\
		\lambda_{jk}x_{jk}^{L} \leq \nu_{jk} \leq \lambda_{jk}x_{jk}^{U}, \;\; j \in D_{k}, \; k \in K \\
		\displaystyle{\sum_{j \in D_{k}} \lambda_{jk} = 1}, \;\; k \in K \\
		\lambda_{jk}h_{jk}(\nu_{jk}/\lambda_{jk}) \leq 0, \;\; j \in D_{k}, \; k \in K \\
		A\lambda \leq a \\
		\lambda_{jk} \in \{0,1\}, \;\; j \in D_{k}, \; k \in K
	\end{array}
\eeq
In this, the Boolean variables $Y_{jk}$ are replaced with binary variables $\lambda_{jk}$. Also, $A\lambda \leq a$ are the logic propositions $\Omega(Y)$ expressed as inequalities in terms of $\lambda$. Notice that $(\lambda_{jk} = 0) \Rightarrow (\nu_{jk} = 0)$, and thus the $j$th system of inequalities in the $k$th disjunction is redundant.

For implementation purposes it is necessary to reformulate (\ref{eq:GDPCHreform}) in such a way so as to avoid division by zero in the nonlinear inequalities. Grossmann and Lee \cite{gross2003} proposed to approximate the set of constraints $\lambda_{jk}h_{jk}(\nu_{jk}/\lambda_{jk}) \leq 0$ by $(\lambda_{jk} + \epsilon)h_{jk}(\nu_{jk}/(\lambda_{jk} + \epsilon)) \leq 0$, where $\epsilon$ is small tolerance, and proved that the approximated constraints are continuous and differentiable. While this transformation is exact for the limiting case when $\epsilon$ tends to zero, the resulting problem risks not being equivalent to the original problem in the sense that the latter's optimal solution would not correspond to the former. In the original problem, when optimal $\lambda_{jk}$ is zero (i.e. $\lambda_{jk}^{*} = 0$), then $\lambda_{jk}^{*}h_{jk}(\nu_{jk}^{*}/\lambda_{jk}^{*})$ becomes zero as well. However, the additional term $\epsilon h_{jk}(\nu_{jk}/(\lambda_{jk} + \epsilon))$ in the approximating constraint prevents it from being equal to 0 when $\lambda_{jk}^{*} = 0$ since $(\lambda_{jk} + \epsilon)h_{jk}(\nu_{jk}/(\lambda_{jk} + \epsilon)) = (0 + \epsilon)h_{jk}(0/(0 + \epsilon)) = \epsilon h_{jk}(0) \neq 0$ since $h_{jk}(0)$ need not be zero. In order to circumvent the feasibility problem, one could attempt to reduce $\epsilon$ to a value small enough such that $\epsilon h_{jk}(0) \leq \tau$ in order to numerically satisfy the constraint within the solver tolerance $\tau$. But this can lead to numerical difficulties since it is not uncommon to require values of $\epsilon$ to be of the order of $10^{-15}$ in order to maintain feasibility. We also note that an approximation of the form $\lambda_{jk}h_{jk}(\nu_{jk}/(\lambda_{jk} + \epsilon)) \leq 0$ does avoid the division by zero and infeasibility problems, but transforms the original convex constraints into nonconvex constraints, which may lead to sub-optimal solutions.

In order to circumvent the issues described above with the approximation proposed by Grossmann and Lee, the following two modifications were proposed by Sawaya and Grossmann \cite{sawaya2007} for the inequality constraints $\lambda_{jk}h_{jk}(\nu_{jk}/\lambda_{jk}) \leq 0$:
\ben
	\item $(\lambda_{jk} + \epsilon)h_{jk}(\nu_{jk}/(\lambda_{jk} + \epsilon)) - \max\limits_{\nu_{jk},\lambda_{jk}} (\epsilon h_{jk}(\nu_{jk}/(\lambda_{jk} + \epsilon)))$
	\item $(\lambda_{jk} + \epsilon)h_{jk}(\nu_{jk}/(\lambda_{jk} + \epsilon)) + h_{jk}(0)(\lambda_{jk} - 1)$
\een
Both modifications resolve the three issues present in the approximation proposed by Lee and Grossmann. In other words,
\bit
	\item they are exact approximations of the original set of constraints $\lambda_{jk}h_{jk}(\nu_{jk}/\lambda_{jk}) \leq 0$ at $\lambda_{jk} = 0$ or $1$ as $\epsilon \rightarrow 0$, 
	\item they maintain feasibility of constraints $\lambda_{jk}h_{jk}(\nu_{jk}/\lambda_{jk}) \leq 0$ at $\lambda_{jk}^{*} = 0$, and 
	\item both approximating functions maintain convexity.
\eit
Although Sawaya's approximating functions address these issues, they are still plagued by the following two drawbacks:
\bit
	\item The second term in the first approximation contains a $max$ function that makes it non-differentiable. For this, we either need extra binary variables or need to solve a global maximization problem, which renders it difficult to implement.
	\item Both approximating functions still require choosing a tolerance level $\epsilon$. The optimal solution can vary substantially with the level of tolerance. Sawaya and Grossmann \cite{sawaya2007} demonstrated fairly different optimal solutions for different values of $\epsilon$ in their work.
\eit

In this work, we propose a novel reformulation methodology for general equality and inequality constrained GDP problems of the form (\ref{eq:GDPgeneral}). The paper is organized as follows. The next section introduces the new methodology with an example. In Section 3, we present its extensions for different kinds of disjunctions derived from implication logic. Section 4 describes its practical applicability. Finally, we present a summary.

\section{New GDP Reformulation Methodology}

The new proposed methodology for reformulating GDP problems into MINLP is applicable for all generic disjunctions represented in the following manner:

\beq
\label{eq:GDPnewmethod}
	\begin{array}{l}
		\bigvee\limits_{j \in D_{k}}
		\left [
			\begin{array}{c}
				Y_{jk} \\
				h_{jk}(x) \leq 0
			\end{array}
		\right ], \quad k \in K \\
		\Omega(Y) = True \\
		x^{L} \leq x \leq x^{U}, \;\;\; Y_{jk} \in \{True,False\} 
	\end{array}
\eeq
Here $x$ is a set of variables that depend on which term in a disjunction is true, such as the mechanical property $P$ above. The constraint set $h_{jk}(x) \leq 0$ includes both equality and inequality constraints. We propose the following new reformulation methodology for such disjunctions:

\beq
\label{eq:newreformulation}
	\begin{array}{l}
		\hat{x}_{jk} = \nu_{jk}^{t} + \nu_{jk}^{f}, \;\; j \in D_{k}, \; k \in K \\
		\lambda_{jk}x^{L} \leq \nu_{jk}^{t} \leq \lambda_{jk}x^{U}, \;\; j \in D_{k}, \; k \in K \\
		(1 - \lambda_{jk})x^{L} \leq \nu_{jk}^{f} \leq (1 - \lambda_{jk})x^{U}, \;\; j \in D_{k}, \; k \in K \\
		\displaystyle{\sum_{j \in D_{k}} \lambda_{jk} = 1}, \;\; k \in K \\
		 h_{jk}(\hat{x}_{jk}) \leq 0, \;\; j \in D_{k}, \; k \in K \\
		\displaystyle{x = \sum_{j \in D_{k}} \nu_{jk}^{t}}, \;\; k \in K \\
		A\lambda \leq a \\
		\lambda_{jk} \in \{0,1\}, \;\; j \in D_{k}, \; k \in K
	\end{array}
\eeq
Here we define artificial variables $\hat{x}_{jk}$ for the set of variables $x$ which are then used in the equation $h_{jk}(\hat{x}_{jk}) \leq 0$. Each $\hat{x}_{jk}$ is disaggregated in a ``true'' variable $\nu_{jk}^{t}$ and a ``false'' variable $\nu_{jk}^{f}$. If a term in a disjunction $k$ is true (i.e. $\lambda_{jk} = 1$), corresponding ``true'' variable $\nu_{jk}^{t}$ is active while the ``false'' variable $\nu_{jk}^{f}$ is set to zero. Similarly, if a term in a disjunction $k$ is not true, the ``false'' variable $\nu_{jk}^{f}$ is active while the ``true'' variable $\nu_{jk}^{t}$ is set to zero. Since only one term in a disjunction $k$ can be true, only one $\nu_{jk}^{t}$ is active while all other $\nu_{jk}^{t}$ are zero. Thus, the value of $\hat{x}_{jk}$ gets assigned to the active $\nu_{jk}^{t}$ which eventually gets transferred to the variable $x$ through the equation $x = \sum_{j \in D_{k}} \nu_{jk}^{t}$. We call this reformulation the ``True-False Reformulation''.

To illustrate with an example, consider the following single disjunction with three terms for each of the inequalities $h_{1}(x)$, $h_{2}(x)$, and $h_{3}(x)$ defined for the variable $x$:

\begin{equation*}
	\begin{array}{l}
		\left [ 
			\begin{array}{c}
				Y_{1} \\
				h_{1}(x) \leq 0
			\end{array}
		\right ] \vee 
		\left [ 
			\begin{array}{c}
				Y_{2} \\
				h_{2}(x) \leq 0
			\end{array}
		\right ] \vee
		\left [ 
			\begin{array}{c}
				Y_{3} \\
				h_{3}(x) \leq 0
			\end{array}
		\right ] \\
		0 \leq x \leq x^{U}
	\end{array}
\end{equation*}
According to the methodology proposed by Lee and Grossmann, this disjunction can be reformulated as follows:
\begin{equation*}
	\begin{array}{l}
		x = \nu_{1} + \nu_{2} + \nu_{3} \\
		0 \leq \nu_{j} \leq \lambda_{j}x^{U}, \;\; j \in \{1,2,3\} \\
		\lambda_{1} + \lambda_{2} + \lambda_{3} = 1 \\
		(\lambda_{1} + \epsilon)h_{1}(\nu_{1}/(\lambda_{1} + \epsilon)) \leq 0 \\
		(\lambda_{2} + \epsilon)h_{2}(\nu_{2}/(\lambda_{2} + \epsilon)) \leq 0 \\
		(\lambda_{3} + \epsilon)h_{3}(\nu_{3}/(\lambda_{3} + \epsilon)) \leq 0 \\
		\lambda_{j} \in \{0,1\}, \;\; j \in \{1,2,3\}
	\end{array}
\end{equation*}
Using the new methodology, we reformulate the disjunction in the following manner:

\begin{equation*}
	\begin{array}{l}
		\hat{x}_{1} = \nu_{1}^{t} + \nu_{1}^{f}, \;\;\;
		\hat{x}_{2} = \nu_{2}^{t} + \nu_{2}^{f}, \;\;\;
		\hat{x}_{3} = \nu_{3}^{t} + \nu_{3}^{f} \\
		0 \leq \nu_{j}^{t} \leq \lambda_{j}x^{U}, \;\; j \in \{1,2,3\} \\
		0 \leq \nu_{j}^{f} \leq (1-\lambda_{j})x^{U}, \;\; j \in \{1,2,3\} \\
		\lambda_{1} + \lambda_{2} + \lambda_{3} = 1 \\
		h_{j}(\hat{x}_{j}) \leq 0, \;\; j \in \{1,2,3\} \\
		x = \nu_{1}^{t} + \nu_{2}^{t} + \nu_{3}^{t} \\
		\lambda_{j} \in \{0,1\}, \;\; j \in \{1,2,3\}
	\end{array}
\end{equation*}
Here we define $\hat{x}_{1}$, $\hat{x}_{2}$, and $\hat{x}_{3}$ as copies of the variable $x$ for each term in the disjunction. If the first term in the disjunction is true, $\lambda_1 = 1$, and we get $\nu_{2}^{t} = 0$ and $\nu_{3}^{t} = 0$. Thus, $x = \nu_{1}^{t}$, which implies $x = \hat{x}_{1}$. Consequently, $h_{1}(\hat{x}_{1}) \leq 0$ ensures that $h_{1}(x) \leq 0$ in the first term in the disjunction is satisfied. The other two inequalities, although still present in the model, do not affect the optimal solution since they are written for dummy ``false'' variables, $h_{2}(\nu_{2}^{f}) \leq 0$ and $h_{3}(\nu_{3}^{f}) \leq 0$.

Compared to a Big-M reformulation, the convex hull reformulation of Lee and Grossmann adds another $n \times \sum_{k} m_{k}$ variables and $n \times q + n \times \sum_{k} m_{k}$ constraints to a GDP problem, where $n$ is the dimension of vector $p$, $m_{k}$ is the number of terms in $k^{th}$ disjunction ($m_{k} = |D_{k}|$), and $q$ is number of disjunctions ($q = |K|$). In comparison, the True-False Reformulation adds $3 \times n \times \sum_{k} m_{k}$ variables and $n \times q + 3 \times n \times \sum_{k} m_{k}$ constraints to the optimization problem. Although we introduce more variables and constraints compared to other reformulations, our proposed methodology provides an improvement over reformulation strategies proposed in the literature since it enjoys all of the following properties, one or more of which are not satisfied by the convex hull or Big-M reformulation methods:
\bit
	\item It provides an exact reformulation of the original disjunctions.
	\item Feasibility of the constraints is maintained, especially at $\lambda_{jk}^{*} = 0$.
	\item Convexity of the constraints is maintained.
	\item The reformulation does not introduce any non-differentiable terms.
	\item The reformulation does not require choosing a tolerance level $\epsilon$.
	\item The reformulation does not require choosing a Big-M parameter.
\eit

We note that the True-False Reformulation methodology is not a replacement for the strategies proposed in the literature. We neither claim nor intend to prove that the new method is relaxed or tighter compared to the Big-M method or Lee's convex hull reformulation approach. Also, we do not prove if the proposed disaggregation of variables in the True-False Reformulation results in a convex hull of the disjunctions. Our intent is just to propose an alternative methodology to convert disjunctions into a standard MINLP format.

\section{Implication Logic Conversion using the New GDP Reformulation}

In this section we demonstrate how the True-False Reformulation can be used to model if-else implication logic by first converting them into mathematical disjunctions and eventually reformulating them to MINLP problems.

\subsection{Converting Implication to Standard MINLP}

Consider the following simplest (but generic) if-else statement that gives rise to conditional constraints. Here either $p = f_{1}(z)$ or $p = f_{2}(z)$ depending on the value taken by the function $g(x)$.\\
\hspace*{1cm}if $g(x) \leq 0$ then  $p = f_{1}(z)$\\
\hspace*{1cm}else $p = f_{2}(z)$\\
This if-else statement can be written in the following implication form:\\
$g(x) \leq 0 \Rightarrow p = f_{1}(z)$\\
$g(x) \nleq 0 \Rightarrow p = f_{2}(z)$\\
This implication can then be transformed to the following logic statements using a negation operator ($\urcorner$):\\
$\urcorner g(x) \leq 0 \vee p = f_{1}(z)$ same as $g(x) \geq 0 \vee p = f_{1}(z)$\\
$\urcorner g(x) \nleq 0 \vee p = f_{2}(z)$ same as $g(x) \leq 0 \vee p = f_{2}(z)$\\
Vecchietti and Grossmann \cite{vecchietti2000} showed that these logic statements can be combined and converted into the following disjunction with two terms:

\begin{equation*}
	\left [
		\begin{array}{c}
			Y_{1} \\
			g(x) \geq 0 \\
			p = f_{2}(z)
		\end{array}
	\right ] \bigvee
	\left [
		\begin{array}{c}
			Y_{2} \\
			g(x) \leq 0 \\
			p = f_{1}(z)
		\end{array}
	\right ]
\end{equation*}
Note that both the testing condition ($g(x) \leq 0$) and the statements in if-else blocks are combined into a single disjunction. Finally, this disjunction can be converted to a standard MINLP format.  Using the True-False Reformulation, it can be transformed in the following manner. For simplicity, we assume the variables $z$ in $f_{1}(z)$ and $f_{2}(z)$ do not overlap with the variables $x$ in $g(x)$.
\begin{equation*}
	\begin{array}{l l}
		\tm{For } x, g(x)
		& \hat{x}_{1} = \nu_{1}^{t} + \nu_{1}^{f}, \;\;\;
		\hat{x}_{2} = \nu_{2}^{t} + \nu_{2}^{f} \\
		& \lambda_{1}x^{L} \leq \nu_{1}^{t} \leq \lambda_{1}x^{U}, \;\; 
		\lambda_{2}x^{L} \leq \nu_{2}^{t} \leq \lambda_{2}x^{U}, \\
		& (1-\lambda_{1})x^{L} \leq \nu_{1}^{f} \leq (1-\lambda_{1})x^{U}, \;\;
		(1-\lambda_{2})x^{L} \leq \nu_{2}^{f} \leq (1-\lambda_{2})x^{U} \\
		& g(\hat{x}_{1}) \geq 0, \;\;\; g(\hat{x}_{2}) \leq 0 \\
		& x = \nu_{1}^{t} + \nu_{2}^{t} \\
		\tm{For } p, f_{1}(z), f_{2}(z)
		& \hat{p}_{1} = \nu_{3}^{t} + \nu_{3}^{f}, \;\;\;
		\hat{p}_{2} = \nu_{4}^{t} + \nu_{4}^{f} \\
		& \lambda_{1}P^{L} \leq \nu_{3}^{t} \leq \lambda_{1}P^{U}, \;\; 
		\lambda_{2}P^{L} \leq \nu_{4}^{t} \leq \lambda_{2}P^{U}, \\
		& (1-\lambda_{1})P^{L} \leq \nu_{3}^{f} \leq (1-\lambda_{1})P^{U}, \;\;
		(1-\lambda_{2})P^{L} \leq \nu_{4}^{f} \leq (1-\lambda_{2})P^{U} \\
		& \hat{p}_{1} = f_{2}(z), \;\;\; \hat{p}_{2} = f_{1}(z) \\
		& p = \nu_{3}^{t} + \nu_{4}^{t} \\
		\lambda_{1} + \lambda_{2} = 1 \\
		\lambda_{j} \in \{0,1\}, \;\; j \in \{1,2\}
	\end{array}
\end{equation*}
The methodology to convert implication construct into standard MINLP problem statements can be summarized in the following two steps:
\bit
	\item \textbf{Step 1:} Define a Boolean variable for each if, else-if, and else sub-blocks in the if-else block. Write each if, else-if, and else sub-block as terms of a single disjunction. Combine the testing condition and block statements in the same term of the disjunction.
	\item \textbf{Step 2:} Reformulate the single disjunction with the True-False Reformulation, disaggregating variables that depend on the if-else construct.
\eit
We illustrate these steps with the help of the following simple example.

\noindent \textbf{Example:} Given the following empirical correlations for Energy consumption ($E$) and Power Cost ($PC$), convert them into a standard MINLP format.\\
\hspace*{1cm}$E = rx + \alpha$\\
\hspace*{1cm}if $E \geq \alpha$ then $PC = PC_{0} + ax^{2} + bx + E - \alpha$\\
\hspace*{1cm}else if $E \leq \beta$ then $PC = PC_{0} - m(\beta - E)$\\
\hspace*{1cm}else $PC = PC_{0}$

\noindent \textbf{Step 1:} We define a Boolean variable $Y_1$, $Y_2$, and $Y_3$ and binary variables $\lambda_1$, $\lambda_2$, and $\lambda_3$ for if, else if, and else blocks, respectively. Next, we write each if, else-if, and else sub-block as terms of a single disjunction.
\begin{equation*}
	\left [
		\begin{array}{c}
			Y_{1} \\
			E \geq \alpha \\
			PC = PC_{0} + ax^{2} + bx + E - \alpha
		\end{array}
	\right ] \vee
	\left [
		\begin{array}{c}
			Y_{2} \\
			E \leq \beta \\
			PC = PC_{0} - m(\beta - E)
		\end{array}
	\right ] \vee
	\left [
		\begin{array}{c}
			Y_{3} \\
			\beta \leq E \leq \alpha \\
			PC = PC_{0}
		\end{array}
	\right ]
\end{equation*}

\noindent \textbf{Step 2:} Finally, we disaggregate variables $E$ and $PC$ (and not the variable $x$), and obtain the following statements using the True-False Reformulation.
\begin{equation*}
	\begin{array}{l l}
		\tm{For $E$: }
		& E = rx + \alpha \\
		& \hat{y}_{j} = \nu_{j}^{t} + \nu_{j}^{f} \;\; j \in \{1,2,3\} \\
		& \hat{y}_{1} \geq \alpha \\
		& \hat{y}_{2} \leq \beta \\
		& \beta \leq \hat{y}_{3} \leq \alpha \\
		& \lambda_{j}E^{L} \leq \nu_{j}^{t} \leq \lambda_{j}E^{U} \;\; j \in \{1,2,3\} \\
		& (1-\lambda_{j})E^{L} \leq \nu_{j}^{f} \leq (1-\lambda_{j})E^{U} \;\; j \in \{1,2,3\} \\
		& \lambda_{1} + \lambda_{2} + \lambda_{3} = 1 \\
		& E = \nu_{1}^{t} + \nu_{2}^{t} + \nu_{3}^{t} \\
		& \lambda_{j} \in \{0,1\}, \;\; j \in \{1,2,3\}
	\end{array}
\end{equation*}

\begin{equation*}
	\begin{array}{l l}
		\tm{For $PC$: }
		& \hat{t}_{j} = \nu_{j}^{t} + \nu_{j}^{f} \;\; j \in \{4,5,6\} \\
		& \hat{t}_{4} = PC_{0} + ax^{2} + bx + \hat{y}_{1} - \alpha \\
		& \hat{t}_{5} = PC_{0} - m(\beta - \hat{y}_{2}) \\
		& \hat{t}_{6} = PC_{0} \\				
		& \lambda_{j}PC^{L} \leq \nu_{j+3}^{t} \leq \lambda_{j}PC^{U} \;\; j \in \{1,2,3\} \\
		& (1-\lambda_{j})PC^{L} \leq \nu_{j+3}^{f} \leq (1-\lambda_{j})PC^{U} \;\; j \in \{1,2,3\} \\
		& PC = \nu_{4}^{t} + \nu_{5}^{t} + \nu_{6}^{t} \\
	\end{array}
\end{equation*}

\subsection{Variations of Implication Logic}

We now discuss three different variations of the implication logic and how they can be tailored to enable application of the True-False Reformulation.

\subsubsection{Implicit else block and sequential implication}

Many implication constructs do not explicitly include an else block. This is usually true when a default case exists irrespective of the implication construct. For instance, the following if-statement does not include an else block because a default equation $g(x)$ exists for the variable $p$.\\
\hspace*{1cm}$p = g(x)$\\
\hspace*{1cm}if $p \leq \alpha$ then  $p = p + f(x)$\\
Here, if $p \leq \alpha$, its value increases by $f(x)$, otherwise remains equal to $g(x)$. In order to convert such a set of statements into disjunctions and standard MINLP format, we introduce dummy variables and an artificial else loop. In particular, the set of statements above can be converted to a standard implication construct with the help of an artificial variable $p_{dummy}$.\\
\hspace*{1cm}$p_{dummy} = g(x)$\\
\hspace*{1cm}if $p_{dummy} \leq \alpha$ then  $p = p_{dummy} + f(x)$\\
\hspace*{1cm}else $p = p_{dummy}$\\
This can then be reformulated into disjunctions and MINLP statements as below:

\begin{equation*}
	\begin{array}{l}
		p_{dummy} = g(x) \\
		\left [
			\begin{array}{c}
				Y_{1} \\
				p_{dummy} \leq \alpha \\
				p = p_{dummy} + f(x)
			\end{array}
		\right ] \vee
		\left [
			\begin{array}{c}
				Y_{2} \\
				p_{dummy} \geq \alpha \\
				p = p_{dummy}
			\end{array}
		\right ] \\
	\end{array}
\end{equation*}

\begin{equation*}
	\begin{array}{l l}
		\tm{For $p_{dummy}$: }
		& p_{dummy} = g(x) \\
		& \hat{p}_{d1} = \nu_{1}^{t} + \nu_{1}^{f}, \;\;\;
		\hat{p}_{d2} = \nu_{2}^{t} + \nu_{2}^{f} \\
		& \hat{p}_{d1} \leq \alpha, \;\;\; \hat{p}_{d2} \geq \alpha \\
		& \lambda P^{L} \leq \nu_{1}^{t} \leq \lambda P^{U}, \;\; 
		(1 - \lambda)P^{L} \leq \nu_{2}^{t} \leq (1 - \lambda)P^{U}, \\
		& (1-\lambda)P^{L} \leq \nu_{1}^{f} \leq (1-\lambda)P^{U}, \;\;
		\lambda P^{L} \leq \nu_{2}^{f} \leq \lambda P^{U} \\
		& p_{dummy} = \nu_{1}^{t} + \nu_{2}^{t} \\
		\tm{For $p$: }
		& \hat{p}_{1} = \nu_{3}^{t} + \nu_{3}^{f}, \;\;\;
		\hat{p}_{2} = \nu_{4}^{t} + \nu_{4}^{f} \\
		& \hat{p}_{1} = \hat{p}_{d1} + f(x), \;\;\; \hat{p}_{2} = \hat{p}_{d2} \\
		& \lambda P^{L} \leq \nu_{3}^{t} \leq \lambda P^{U}, \;\; 
		(1 - \lambda)P^{L} \leq \nu_{4}^{t} \leq (1 - \lambda)P^{U}, \\
		& (1-\lambda)P^{L} \leq \nu_{3}^{f} \leq (1-\lambda)P^{U}, \;\;
		\lambda P^{L} \leq \nu_{4}^{f} \leq \lambda P^{U} \\
		& p = \nu_{3}^{t} + \nu_{4}^{t} \\
	\end{array}
\end{equation*}

Another related case is of sequential statements. In other words, the value of a variable can depend on multiple implication constructs running in sequence. For example, in the following, the value of $p$ obtained is eventually decided after three sequential if-else blocks.\\
\hspace*{1cm}if $r \leq \alpha$ then $p = f_{1}(x)$\\
\hspace*{1cm}else $p = f_{2}(x)$\\
\hspace*{1cm}if $p \geq \beta$ then $p = \beta$\\
\hspace*{1cm}if $p \leq \gamma$ then $p = \gamma$\\
For multiple sequential implication blocks, we define multiple artificial variables and else blocks. For the example above, we introduce two additional dummy variables $p_{d}^{1}$ and $p_{d}^{2}$. Using these, the implication logic is then converted to the following with dummy else blocks\\
\hspace*{1cm}if $r \leq \alpha$ then $p_{d}^{1} = f_{1}(x)$\\
\hspace*{1cm}else $p_{d}^{1} = f_{2}(x)$\\
\hspace*{1cm}if $p_{d}^{1} \geq \beta$ then $p_{d}^{2} = \beta$\\
\hspace*{1cm}else $p_{d}^{2} = p_{d}^{1}$\\
\hspace*{1cm}if $p_{d}^{2} \leq \gamma$ then $p = \gamma$\\
\hspace*{1cm}else $p = p_{d}^{2}$\\
Finally, this can be reformulated into the following three disjunctions

\begin{equation*}
	\left [
		\begin{array}{c}
			Y_{1} \\
			r \leq \alpha \\
			p_{d}^{1} = f_{1}(x)
		\end{array}
	\right ] \vee
	\left [
		\begin{array}{c}
			Y_{2} \\
			r \geq \alpha \\
			p_{d}^{1} = f_{2}(x)
		\end{array}
	\right ]
\end{equation*}
\begin{equation*}
	\left [
		\begin{array}{c}
			Z_{1} \\
			p_{d}^{1} \geq \beta \\
			p_{d}^{2} = \beta
		\end{array}
	\right ] \vee
	\left [
		\begin{array}{c}
			Z_{2} \\
			p_{d}^{1} \leq \beta \\
			p_{d}^{2} = p_{d}^{1}
		\end{array}
	\right ]
\end{equation*}
\begin{equation*}
	\left [
		\begin{array}{c}
			W_{1} \\
			p_{d}^{2} \leq \gamma \\
			p = \gamma
		\end{array}
	\right ] \vee
	\left [
		\begin{array}{c}
			W_{2} \\
			p_{d}^{2} \geq \gamma \\
			p = p_{d}^{2}
		\end{array}
	\right ]
\end{equation*}
These disjunctions can be converted to the equations of an MINLP problem using the steps of the True-False Reformulation explained above.

\subsubsection{Multiple testing conditions}

The cases considered so far included only a single testing condition in the implication logic. On many occasions we encounter multiple testing conditions in the same implication construct. For example, consider the following set of statements which includes two testing conditions $E \geq \beta$ and $E \leq \alpha$ in the same if-else block:\\
\hspace*{1cm}if $E \geq \beta$ and $E \leq \alpha$ then  $PC = f_{1}(x)$\\
\hspace*{1cm}else $PC = f_{2}(x)$\\
In order to apply the True-False Reformulation methodology, we modify the Step 1 of the algorithm for such cases. In particular, we make the following two changes to Step 1:
\ben
	\item In addition to defining a binary variable (and Boolean variable) for each if, else if, and else blocks, we also define binary (and Boolean) variables for each testing condition.
	\item We define separate disjunctions for the testing conditions and the statements in each block instead of combining them in a single disjunction.
\een
In Step 2, the testing conditions are reformulated using the binary variables defined for testing conditions, while the inequalities in the if-else blocks are reformulated using the block binary variables. The two sets of binary variables are linked via implications and logical constraints. For example, the aforementioned if-else statements can be converted into the following disjunctions. Here we define Boolean variables $Z_{1}$ and $Z_{2}$ for $E \geq \beta$ and $E \leq \alpha$, respectively, while $Y_{1}$ and $Y_{2}$ for the if and else blocks. Both sets of binary variables are linked using implications.

\begin{equation*}
	\begin{array}{l}
		\left [
			\begin{array}{c}
				Z_{1} \\
				E \geq \beta
			\end{array}
		\right ] \vee
		\left [
			\begin{array}{c}
				Z_{2} \\
				E \leq \alpha
			\end{array}
		\right ] \\
		\left [
			\begin{array}{c}
				Y_{1} \\
				PC = f_{1}(x)
			\end{array}
		\right ] \vee
		\left [
			\begin{array}{c}
				Y_{2} \\
				PC = f_{2}(x)
			\end{array}
		\right ] \\
		Z_{1} \wedge Z_{2} \Rightarrow Y_{1} \tm{ which is same as } 
		\urcorner Z_{1} \vee \urcorner Z_{2} \vee Y_{1} \\
		\urcorner(Z_{1} \wedge Z_{2}) \Rightarrow Y_{2} \tm{ which is same as } 
			\left \{
				\begin{array}{l}
					Z_{1} \vee Y_{2} \\
					Z_{2} \vee Y_{2}
				\end{array}
			\right .			
	\end{array}
\end{equation*}
Here ($\wedge$) is the logical ``and'' operator. The disjunctions are reformulated in the following MINLP statements using the True-False Reformulation. Here, we reformulate the testing conditions $E \geq \beta$ and $E \leq \alpha$ using their binary variables $z_1$ and $z_2$ (for Boolean variables $Z_1$ and $Z_2$, respectively). We convert implications into mathematical constraints to connect the binary variable $\lambda_1$ and $\lambda_2$ (for $Y_1$ and $Y_2$, respectively) to $z_1$ and $z_2$.
\begin{equation*}
	\begin{array}{l l}
		\tm{Testing conditions: }
		& \hat{y}_{1} = \nu_{1}^{t} + \nu_{1}^{f}, \;\;\;
		\hat{y}_{2} = \nu_{2}^{t} + \nu_{2}^{f} \\
		& \hat{y}_{1} \geq \beta, \;\;\; \hat{y}_{2} \leq \alpha \\
		& z_{1}E^{L} \leq \nu_{1}^{t} \leq z_{1}E^{U}, \;\;
		z_{2}E^{L} \leq \nu_{2}^{t} \leq z_{2}E^{U}, \\
		& (1-z_{1})E^{L} \leq \nu_{1}^{f} \leq (1-z_{1})E^{U}, \;\;
		(1-z_{2})E^{L} \leq \nu_{2}^{f} \leq (1-z_{2})E^{U} \\
		& z_{1} + z_{2} = 1 \\
		& E = \nu_{1}^{t} + \nu_{2}^{t} \\
		\tm{If-else blocks: }
		& \hat{p}_{1} = \nu_{3}^{t} + \nu_{3}^{f}, \;\;\;
		\hat{p}_{2} = \nu_{4}^{t} + \nu_{4}^{f} \\
		& \hat{p}_{1} = f_{1}(x), \;\;\; \hat{p}_{2} = f_{2}(x) \\
		& \lambda_{1}PC^{L} \leq \nu_{3}^{t} \leq \lambda_{1}PC^{U}, \;\;
		\lambda_{2}PC^{L} \leq \nu_{4}^{t} \leq \lambda_{2}PC^{U}, \\
		& (1-\lambda_{1})PC^{L} \leq \nu_{3}^{f} \leq (1-\lambda_{1})PC^{U}, \;\;
		(1-\lambda_{2})PC^{L} \leq \nu_{4}^{f} \leq (1-\lambda_{2})PC^{U} \\
		& \lambda_{1} + \lambda_{2} = 1 \\
		& PC = \nu_{3}^{t} + \nu_{4}^{t} \\
		\tm{Implications: }
		& \lambda_{1} \geq z_{1} + z_{2} - 1, \;\;\; 
		\lambda_{2} + z_{1} \geq 1, \;\;\;
		\lambda_{2} + z_{2} \geq 1 \\
		& \lambda_{1}, \lambda_{2}, z_{1}, z_{2} \in \{0,1\}
	\end{array}
\end{equation*}

We note that this framework can be used to break down any complex arrangement of testing conditions in implications. For instance, consider the following if-else construct\\
\hspace*{1cm}if ($p_{1} \leq \alpha$ and $p_{2} \leq \beta$) or ($p_{1} \geq \gamma$ and $p_{2} \geq \delta$) then \\ 
\hspace*{2cm}$PC = f_{1}(x)$\\
\hspace*{1cm}else\\
\hspace*{2cm}$PC = f_{2}(x)$\\
\hspace*{1cm}end\\
This can be reformulated into the following set of disjunctions and logical statements.

\begin{equation*}
	\begin{array}{l}
		\left [
			\begin{array}{c}
				Z_{1} \\
				p_{1} \leq \alpha
			\end{array}
		\right ] \vee
		\left [
			\begin{array}{c}
				Z_{2} \\
				p_{2} \leq \beta
			\end{array}
		\right ] \vee
		\left [
			\begin{array}{c}
				Z_{3} \\
				p_{1} \geq \gamma
			\end{array}
		\right ] \vee
		\left [
			\begin{array}{c}
				Z_{4} \\
				p_{2} \geq \delta
			\end{array}
		\right ] \\
		\left [
			\begin{array}{c}
				Y_{1} \\
				PC = f_{1}(x)
			\end{array}
		\right ] \vee
		\left [
			\begin{array}{c}
				Y_{2} \\
				PC = f_{2}(x)
			\end{array}
		\right ] \\
		(Z_{1} \wedge Z_{2}) \vee (Z_{3} \wedge Z_{4}) \Rightarrow Y_{1} 
		\tm{ which is same as }
		\left \{
			\begin{array}{l}
				\urcorner Z_{1} \vee \urcorner Z_{2} \vee Y_{1} \\
				\urcorner Z_{3} \vee \urcorner Z_{4} \vee Y_{1}
			\end{array}
		\right . \\
		\urcorner((Z_{1} \wedge Z_{2}) \vee (Z_{3} \wedge Z_{4})) \Rightarrow Y_{2} 
		\tm{ which is same as } 
		\left \{
			\begin{array}{l}
				Z_{1} \vee Z_{3} \vee Y_{2} \\
				Z_{2} \vee Z_{3} \vee Y_{2} \\
				Z_{1} \vee Z_{4} \vee Y_{2} \\
				Z_{2} \vee Z_{4} \vee Y_{2}
			\end{array}
		\right .			
	\end{array}
\end{equation*}

\subsubsection{Nested implications}

Finally we consider the nested implications. In order to convert nested implication logic blocks into disjunctions, we follow a similar treatment used for the sequential implication logic in Section 3.2.1. In particular, we define dummy variables to connect inner and outer if-else blocks. For example, let us consider the following nested if-else statements \\
\hspace*{1cm}if $p_{1} \leq \alpha$ then \\
\hspace*{2cm}$p_{2} = g(x)$ \\
\hspace*{2cm}if $p_{1} \geq \kappa p_{2}$ then $T = f_{1}(x)$\\
\hspace*{2cm}else $T = f_{2}(x)$\\
\hspace*{1cm}else\\
\hspace*{2cm}$T = f_{3}(x)$\\
\hspace*{1cm}end\\
In order to reformulate it into disjunctions, we define an artificial variable $T_d$ for the inner if-else construct. Next, we convert the nested structure into the following sequential structure of implications. Here we first write statements of the inner if-else construct with the dummy variable $T_d$, while the statements of the outer if-else block follow after that.\\
\hspace*{1cm}$p_{2} = g(x)$ \\
\hspace*{1cm}if $p_{1} \geq \kappa p_{2}$ then $T_{d} = f_{1}(x)$\\
\hspace*{1cm}else $T_{d} = f_{2}(x)$\\
\hspace*{1cm}if $p_{1} \leq \alpha$ then $T = T_{d}$\\
\hspace*{1cm}else $T = f_{3}(x)$\\
Note that the latter implication construct is not equivalent to the former. In the former, functions $g(x)$, $f_{1}(x)$, and $f_{2}(x)$ are evaluated only when $p_{1} \leq \alpha$. If $p_{1} > \alpha$, only $f_{3}(x)$ is evaluated. In contrast, in the latter construct, functions $g(x)$, $f_{1}(x)$, and $f_{2}(x)$ and variable $p_2$ are computed irrespective of the value of $p_1$. We can avoid this by choosing an alternate representation where the nested implication blocks are converted into the following if-else-if ladder\\
\hspace*{1cm}$p_{2} = g(x)$ \\
\hspace*{1cm}if $p_{1} \leq \alpha$ and $p_{1} \geq \kappa p_{2}$ then $T = f_{1}(x)$\\
\hspace*{1cm}else if $p_{1} \leq \alpha$ and $p_{1} \leq \kappa p_{2}$ then $T = f_{2}(x)$\\
\hspace*{1cm}else $T = f_{3}(x)$\\
However, this representation requires introducing additional binary variables for multiple testing conditions. Moreover, generating such a ladder can become cumbersome when the nesting increases both in length and depth. Thus, we prefer to use the representation with dummy variables. 

Once converted with the help of dummy variables, the sequential implication construct can then be transformed into the following two disjunctions.

\begin{equation*}
	\begin{array}{l}
		p_{2} = g(x) \\
		\left [
			\begin{array}{c}
				Y_{1} \\
				p_{1} \geq \kappa p_{2} \\
				T_{d} = f_{1}(x)
			\end{array}
		\right ] \vee
		\left [
			\begin{array}{c}
				Y_{2} \\
				p_{1} \leq \kappa p_{2} \\
				T_{d} = f_{2}(x)
			\end{array}
		\right ]
	\end{array}
\end{equation*}
\begin{equation*}
	\begin{array}{l}
		\left [
			\begin{array}{c}
				Z_{1} \\
				p_{1} \leq \alpha \\
				T = T_{d}
			\end{array}
		\right ] \vee
		\left [
			\begin{array}{c}
				Z_{2} \\
				p_{1} \geq \alpha \\
				T = f_{3}(x)
			\end{array}
		\right ]
	\end{array}
\end{equation*}
Finally, the disjunctions can be reformulated to the following equations in the MINLP format using the True-False Reformulation.

\begin{equation*}
	\begin{array}{l l}
		\tm{First disjunction: }
		& p_{2} = g(x) \\
		& \hat{p}_{1a} = \nu_{1a}^{t} + \nu_{1a}^{f}, \;\;\; 
		\hat{p}_{1b} = \nu_{1b}^{t} + \nu_{1b}^{f} \\
		& \hat{p}_{1a} \geq \kappa p_{2}, \;\;\; \hat{p}_{1b} \leq \kappa p_{2} \\
		& y_{1}P^{L} \leq \nu_{1a}^{t} \leq y_{1}P^{U}, \;\;\;
		y_{2}P^{L} \leq \nu_{1b}^{t} \leq y_{2}P^{U}, \\
		& (1-y_{1})P^{L} \leq \nu_{1b}^{f} \leq (1-y_{1})P^{U}, \;\;\;
		(1-y_{2})P^{L} \leq \nu_{1b}^{f} \leq (1-y_{2})P^{U}, \\
		& p_{1} = \nu_{1a}^{t} + \nu_{1b}^{t} \\
		& \hat{t}_{d1} = \nu_{d1}^{t} + \nu_{d1}^{f}, \;\;\; 
		\hat{t}_{d2} = \nu_{d2}^{t} + \nu_{d2}^{f} \\
		& \hat{t}_{d1} = f_{1}(x), \;\;\; \hat{t}_{d2} = f_{2}(x) \\		
		& y_{1}T^{L} \leq \nu_{d1}^{t} \leq y_{1}T^{U}, \;\;\;
		y_{2}T^{L} \leq \nu_{d2}^{t} \leq y_{2}T^{U}, \\
		& (1-y_{1})T^{L} \leq \nu_{d1}^{f} \leq (1-y_{1})T^{U}, \;\;\;
		(1-y_{2})T^{L} \leq \nu_{d2}^{f} \leq (1-y_{2})T^{U}, \\
		& T_{d} = \nu_{d1}^{t} + \nu_{d2}^{t} \\
		& y_{1} + y_{2} = 1 \\
		\tm{Second disjunction: }
		& \hat{p}_{2a} = \nu_{2a}^{t} + \nu_{2a}^{f}, \;\;\; 
		\hat{p}_{2b} = \nu_{2b}^{t} + \nu_{2b}^{f} \\
		& \hat{p}_{2a} \leq \alpha, \;\;\; \hat{p}_{2b} \geq \alpha \\
		& z_{1}P^{L} \leq \nu_{2a}^{t} \leq z_{1}P^{U}, \;\;\;
		z_{2}P^{L} \leq \nu_{2b}^{t} \leq z_{2}P^{U}, \\
		& (1-z_{1})P^{L} \leq \nu_{2b}^{f} \leq (1-z_{1})P^{U}, \;\;\;
		(1-z_{2})P^{L} \leq \nu_{2b}^{f} \leq (1-z_{2})P^{U}, \\
		& p_{1} = \nu_{2a}^{t} + \nu_{2b}^{t} \\
		& \hat{t}_{1} = \nu_{1}^{t} + \nu_{1}^{f}, \;\;\; 
		\hat{t}_{2} = \nu_{2}^{t} + \nu_{2}^{f} \\
		& \hat{t}_{1} = T_{d}, \;\;\; \hat{t}_{2} = f_{3}(x) \\		
		& z_{1}T^{L} \leq \nu_{1}^{t} \leq z_{1}T^{U}, \;\;\;
		z_{2}T^{L} \leq \nu_{2}^{t} \leq z_{2}T^{U}, \\
		& (1-z_{1})T^{L} \leq \nu_{1}^{f} \leq (1-z_{1})T^{U}, \;\;\;
		(1-z_{2})T^{L} \leq \nu_{2}^{f} \leq (1-z_{2})T^{U}, \\
		& T = \nu_{1}^{t} + \nu_{2}^{t} \\
		& z_{1} + z_{2} = 1 \\
		& y_{j}, z_{j} \in \{0,1\} \;\;\; j \in \{1,2\}
	\end{array}
\end{equation*}

\section{Practical Application}

The True-False Reformulation presented in this work has been implemented in a mathematical optimal product design tool that provides a systematic methodology to design commercial formulated products (Section 3 in \cite{agarwal2012}).  For customer-desired targets and restrictions for physical and mechanical properties provided as inputs, this mathematical optimization tool constructs an MINLP problem using nonlinear physical property models, and eventually generates a list of multiple potential products/compounds that satisfy customer property specifications.  

The physical and mechanical property models in the product design tool are defined in terms of multiple implication logic statements.  For example, the fundamental model for one of the mechanical properties $P$ is represented using the following chain of implication statements (note that in the following we do not completely disclose the equations for the mechanical property $P$ for confidentiality reasons)\\
\hspace*{1cm}if ($F_{T} \geq \alpha$) or ($M_{c} \leq \beta$) then\\
\hspace*{2cm}$P = f_{1}(x)$\\
\hspace*{1cm}else\\
\hspace*{2cm}$P = f_{2}(x)$\\
\hspace*{1cm}end\\
\hspace*{1cm}$ebs = g(x)$\\
\hspace*{1cm}if $ebs \geq \gamma$ then\\
\hspace*{2cm}if $Mc \geq 0$ then $ec = h_{1}(x)$\\
\hspace*{2cm}else $ec = h_{2}(x)$\\
\hspace*{2cm}if $ec \geq \rho$ then $ec = \rho$\\
\hspace*{2cm}if $ebs \leq \kappa ec$ then\\
\hspace*{3cm}$P_{corr} = s_{1}(x)$\\
\hspace*{2cm}else if $\kappa ec \leq ebs \leq ec$ then\\
\hspace*{3cm}$P_{corr} = s_{2}(x)$\\
\hspace*{2cm}else\\
\hspace*{3cm}$P_{corr} = s_{3}(x)$\\
\hspace*{2cm}end\\
\hspace*{2cm}if $P_{corr} \leq \delta P$ then $P = \theta P_{corr}$\\
\hspace*{2cm}else $P = \eta (P_{corr} + P)$\\
\hspace*{1cm}end\\
This implication construct comprises all variations of implication logic considered in Section 3; in particular, implicit else blocks, sequential implications, multiple testing conditions, and nested implication blocks.  Consequently, extensions of the True-False Reformulation developed in Section 3 are applied and utilized in the optimal product design tool to reformulate equations of property $P$ in the MINLP format.  We convert the aforementioned implication statements into the following set of disjunctions.

\begin{equation*}
	\begin{array}{l}
		\left [
			\begin{array}{c}
				Z_{1} \\
				F_{T} \geq \alpha
			\end{array}
		\right ] \vee
		\left [
			\begin{array}{c}
				Z_{2} \\
				M_{c} \leq \beta
			\end{array}
		\right ]
	\end{array}
\end{equation*}
\begin{equation*}
	\begin{array}{l}
		\left [
			\begin{array}{c}
				Y_{11} \\
				P_{dummy}^{1} = f_{1}(x)
			\end{array}
		\right ] \vee
		\left [
			\begin{array}{c}
				Y_{12} \\
				P_{dummy}^{1} = f_{2}(x)
			\end{array}
		\right ]
	\end{array}
\end{equation*}
\begin{equation*}
	\begin{array}{l}
		\left [
			\begin{array}{c}
				Y_{21} \\
				Mc \geq 0 \\
				ec_{dummy} = h_{1}(x)
			\end{array}
		\right ] \vee
		\left [
			\begin{array}{c}
				Y_{22} \\
				Mc \leq 0 \\
				ec_{dummy} = h_{2}(x)
			\end{array}
		\right ]
	\end{array}
\end{equation*}
\begin{equation*}
	\begin{array}{l}
		\left [
			\begin{array}{c}
				Y_{31} \\
				ec_{dummy} \geq \rho \\
				ec = \rho
			\end{array}
		\right ] \vee
		\left [
			\begin{array}{c}
				Y_{32} \\
				ec_{dummy} \leq \rho \\
				ec = ec_{dummy}
			\end{array}
		\right ]
	\end{array}
\end{equation*}
\begin{equation*}
	\begin{array}{l}
		\left [
			\begin{array}{c}
				Y_{41} \\
				ebs \leq \kappa ec \\
				P_{corr} = s_{1}(x)
			\end{array}
		\right ] \vee
		\left [
			\begin{array}{c}
				Y_{42} \\
				\kappa ec \leq ebs \leq ec \\
				P_{corr} = s_{2}(x)
			\end{array}
		\right ] \vee
		\left [
			\begin{array}{c}
				Y_{43} \\
				ebs \geq ec \\
				P_{corr} = s_{3}(x)
			\end{array}
		\right ]
	\end{array}
\end{equation*}
\begin{equation*}
	\begin{array}{l}
		\left [
			\begin{array}{c}
				Y_{51} \\
				P_{corr} \leq \delta P_{dummy}^{1} \\
				P_{dummy}^{2} = \theta P_{corr}
			\end{array}
		\right ] \vee
		\left [
			\begin{array}{c}
				Y_{52} \\
				P_{corr} \geq \delta P_{dummy}^{1} \\
				P_{dummy}^{2} = \eta (P_{corr} + P_{dummy}^{1})
			\end{array}
		\right ]
	\end{array}
\end{equation*}
\begin{equation*}
	\begin{array}{l}
		\left [
			\begin{array}{c}
				Y_{61} \\
				ebs \geq \gamma \\
				P = P_{dummy}^{1}
			\end{array}
		\right ] \vee
		\left [
			\begin{array}{c}
				Y_{62} \\
				ebs \leq \gamma \\
				P = P_{dummy}^{2}
			\end{array}
		\right ] \\
		(Z_{1} \vee Z_{2}) \Rightarrow Y_{11} 
		\tm{ which is same as }
		\left \{
			\begin{array}{l}
				\urcorner Z_{1} \vee Y_{11} \\
				\urcorner Z_{2} \vee Y_{11}
			\end{array}
		\right . \\
		\urcorner((Z_{1} \vee Z_{2}) \Rightarrow Y_{12} 
		\tm{ which is same as }
		Z_{1} \vee Z_{2} \vee Y_{12}
	\end{array}
\end{equation*}
These disjunctions can then be reformulated into MINLP statements using the True-False Reformulation. Other property models in the optimal product design tool have been similarly treated using the True-False Reformulation.

\section{Summary}

We summarize the key points of this paper as follows:
\bit
	\item Generalized Disjunctive Programming provides a systematic methodology to model optimization problems involving logic disjunctions, logic propositions, and algebraic equations. In order to take advantage of the existing MINLP solvers, GDP problems can be reformulated as the standard MINLP problems. 
	\item In this work we propose a novel reformulation methodology, called the True-False Reformulation, for generic equality and inequality constrained GDP problems. It is a variant of the convex hull reformulation approach.
	\item The proposed approach involves defining an artificial variable for each term in a disjunction. The artificial variable is then disaggregated into a ``true'' and a ``false variable.
	\item The new methodology provides an exact reformulation, maintains feasibility and convexity of the constraints, and does not require choosing a tolerance level $\epsilon$ and a Big-M parameter unlike the reformulation approaches in the literature.
	\item We develop a systematic methodology to convert implication logic into equations in the standard MINLP format using the True-False Reformulation.
	\item The systematic conversion methodology for implication logic comprises two steps:
	\bit
		\renewcommand{\labelitemii}{$\circ$}
		\item Define a Boolean variable for each if, else-if, and else sub-blocks in the if-else block. Write each if, else-if, and else sub-block as terms of a single disjunction. Combine the testing condition and block statements in the same term of the disjunction.
		\item Reformulate the disjunction with the True-False Reformulation, disaggregating variables that depend on the if-else construct.
	\eit
	\item The two-step approach is extended for variations of implication constructs that include implicit else blocks, sequential implications, multiple testing conditions, and nested implication blocks.
\eit

\bibliographystyle{amsplain}

\providecommand{\bysame}{\leavevmode\hbox to3em{\hrulefill}\thinspace}
\providecommand{\MR}{\relax\ifhmode\unskip\space\fi MR }
\providecommand{\MRhref}[2]{%
  \href{http://www.ams.org/mathscinet-getitem?mr=#1}{#2}
}
\providecommand{\href}[2]{#2}

\end{document}